\documentclass[11pt]{article}
\usepackage[T1]{fontenc}
\usepackage{lmodern}
\usepackage[margin=0.75in]{geometry}
\usepackage{amsmath,amssymb,amsthm,mathtools,microtype}
\usepackage[hidelinks]{hyperref}
\usepackage{graphicx}
\usepackage{placeins}
\newtheorem{theorem}{Theorem}[section]
\newtheorem{lemma}[theorem]{Lemma}
\newtheorem{proposition}[theorem]{Proposition}
\newtheorem{corollary}[theorem]{Corollary}
\theoremstyle{remark}
\newtheorem{remark}[theorem]{Remark}
\newtheorem{definition}[theorem]{Definition}
\numberwithin{equation}{section}
\title{Mortgage Backed Securities with p-Laplacian Diffusion and Random Prepayment}
\author{%
    \qquad Dianchao
Wu$^{\dagger}$ and Rana D. Parshad$^{\dagger}$\\[0.8em]
  {\small $^{\dagger}$Department of Mathematics, Iowa State University}\\
  {\small\texttt{dianchao@iastate.edu}\quad\texttt{rparshad@iastate.edu}}\\
}
\date{November 2026}
\begin{document}

\maketitle
\begin{abstract}
    In this note we consider deterministic and stochastic versions of a PDE pricing model for mortgage backed securities, driven by p-Laplacian motion. We show the existence of weak solutions in the $1<p<2$ case, for sufficiently small initial data. In the stochastic case, we consider a random prepayment function, when $p=2$. Herein we derive results concerning the uniform first order sensitivity and mean squared sensitivity of the solution, to the pricing model.
\end{abstract}
\section{Introduction}

A mortgage-backed security (MBS) is a claim on the cash flows of a pool of
mortgages. In a pass-through security, each investor receives a pro rata share of
the interest and principal paid by the borrowers. When mortgage rates fall,
borrowers tend to exercise their right to prepay. The investor then receives
principal early, earns no further interest on it, and may have to reinvest at a
lower rate \cite{TaskForce2003}. The timing of prepayments therefore affects the
value of an MBS, and prepayment modelling is central to its valuation.

There are two main approaches to modelling prepayment. In empirical models, such as
that of Schwartz and Torous \cite{SchwartzTorous1989}, data are used to estimate
how the prepayment rate depends on explanatory variables. In option-based models,
which go back to Dunn and McConnell \cite{DunnMcConnell1981}, borrowers prepay when
it is optimal for them to do so. Stanton \cite{Stanton1995} extended this approach
to borrowers who face heterogeneous transaction costs and may delay refinancing.
Gabaix, Krishnamurthy and Vigneron \cite{GKV2007} proposed an equilibrium model in
which prepayment risk is priced, and tested it using market data. Building on this
equilibrium model, Papi and Briani \cite{PapiBriani2011} derived a partial
differential equation for the price. Let $t$ be the time to maturity and let
$z\in\mathbb{R}^N$ be the vector of economic factors that affect prepayment, with
drift $\mu(z,t)$ and volatility matrix $\sigma(t)$. If $h(z,t)$ denotes the
remaining principal and $P=u+h$ the value of the security, then $u$ satisfies
\begin{equation}\label{eq:intro-pb}
\begin{cases}
u_t=\dfrac12\operatorname{tr}\bigl(\sigma\sigma^{\top}\nabla^{2}u\bigr)
+\mu\cdot\nabla u-\rho\,\dfrac{|\sigma^{\top}\nabla u|^{2}}{u+h+\xi}
+(\tau-r)h-ru, & (z,t)\in\mathbb{R}^{N}\times(0,T),\\[2mm]
u(z,0)=0, & z\in\mathbb{R}^{N},
\end{cases}
\end{equation}
where $\tau$ is the coupon rate, $r$ is the discount rate, $\rho>0$ is a risk
parameter and $\xi(t)=A_0e^{r(T-t)}$ with $A_0>0$.

Equation \eqref{eq:intro-pb} is semilinear with quadratic growth in the gradient.
It is degenerate when some of the factors have no diffusion, as in the applications
considered in \cite{Papi2003,PapiBriani2011}. Papi \cite{Papi2002} proved a
comparison principle for viscosity solutions of \eqref{eq:intro-pb}, from which
existence and uniqueness follow, and the regularity of the solution was examined in
\cite{Papi2003,PapiBriani2011}. In the nondegenerate case with constant
coefficients, Parshad \cite{Parshad2010} showed that, under further assumptions on
the prepayment function, \eqref{eq:intro-pb} can be reduced to a semilinear
equation without the quadratic term, and proved that the reduced model is well
posed on a bounded domain. Strong solutions of the reduced model were obtained in
\cite{PBBP2011}, where the model was also compared against market data.

In this paper we study two questions about \eqref{eq:intro-pb}. The first concerns
the diffusion term. In \eqref{eq:intro-pb} the diffusive flux of $u$ is linear in
$\nabla u$. We replace this linear relation by a power law under which the flux is
parallel to $\nabla u$ and its magnitude is proportional to $|\nabla u|^{p-1}$. In
the isotropic case $\sigma\sigma^{\top}=2aI$ this amounts to replacing $a\Delta u$
by the $p$-Laplacian $a\nabla\cdot(|\nabla u|^{p-2}\nabla u)$. Power laws of this
type are classical in the theory of non-Newtonian fluids and of filtration in
porous media, and the resulting equations have a well-developed theory
\cite{DiBenedetto1993}. In finance, Tak\'a\v{c} \cite{Takac2023} recently proposed
a Black--Scholes model with the $p$-Laplacian for option pricing. 

The analysis of \eqref{eq:intro-pb} relies on an exact match between its two nonlinear features.
The diffusion is linear and the gradient term is quadratic, so that a change of
unknown removes the latter \cite{Parshad2010}, and the comparison principle of
\cite{Papi2002} rests on a compatibility condition of the same kind. Under a power
law this match is lost. The flux is homogeneous of degree $p-1$ in $\nabla u$ while
the gradient term remains quadratic, and for $p<2$ it grows faster than
$|\nabla u|^p$, the natural growth for the $p$-Laplacian. We ask how much of the
theory survives.We focus on the range $1<p<2$, known as the fast diffusion range \cite{BIV2010}, in which the flux
grows sublinearly in the gradient and the equation is singular where $\nabla u=0$.
The resulting equation \eqref{eq:rp-pricing-pde-plap} is posed on a bounded domain
in one or two space dimensions, with homogeneous Neumann boundary conditions.
Unlike \eqref{eq:intro-pb}, it is not derived from an equilibrium argument. For
$3/2<p<2$ and small bounded initial data, we establish estimates for a regularized
problem that hold uniformly in the regularization parameter (Theorem~\ref{thm:mt11}
and Lemma~\ref{lem:ge11n}). We also obtain weak limits of the regularized
solutions along subsequences (Theorem~\ref{thm:vubd}). The main result on weak solutions of \eqref{eq:rp-pricing-pde-plap} is provided via Theorem \ref{thm:mt1}.

The second question concerns the prepayment intensity. We consider the two-factor
case of \eqref{eq:intro-pb} from \cite[Example~10.26]{PapiBriani2011}, in which
$z=(x,y)$, the financial factor $x$ diffuses, and the cumulative prepayment $y$
grows at a rate $\eta(z,t)$, the prepayment intensity. In applications the
intensity is known only approximately because it is estimated from historical data
\cite{SchwartzTorous1989}. To describe this uncertainty, we multiply a baseline
intensity $\eta^0$ by a random factor $1+\varepsilon\zeta(t)$, where the random
function $\zeta$ is continuous in time and bounded by a constant, and we ask how
the price responds for small $\varepsilon$. Let $u^\varepsilon$ denote the solution
with intensity $\eta^0(1+\varepsilon\zeta)$. Our main result,
Theorem~\ref{theorem:rp-first-order}, states that there is a random field $Y$,
which depends linearly on $\zeta$, such that for all sufficiently small
$\varepsilon$, almost surely,
\begin{equation}\label{eq:intro-main}
\sup_{(z,t)\in\mathbb{R}^2\times[0,T]}
\bigl|u^{\varepsilon}(z,t)-u^{0}(z,t)-\varepsilon Y(z,t)\bigr|
\le C\varepsilon^{2}\|\zeta\|_{\infty}^{2},
\end{equation}
with a deterministic constant $C$. As a consequence
(Corollary~\ref{theorem:rp-mean-variance}), if $\zeta$ is centered then the mean
price does not change at first order in $\varepsilon$, and the variance of the
price is $\varepsilon^{2}\operatorname{Var}(Y)$ up to an error of order
$|\varepsilon|^{3}$. Thus, to leading order, the variance depends on the law of
$\zeta$ only through its covariance function.

The proof is probabilistic. In the two-factor case the coefficient $a=\sigma^2/2$
of the diffusion and the coefficient $k=\rho\sigma^2$ of the quadratic term satisfy
$k=2\rho a$, and a power transformation of $u+h+\xi$, of the kind used in
\cite{Parshad2010}, removes the quadratic term. Since this transformation is
unavailable for $p\neq2$, we treat the second question only for linear diffusion.
The transformed equation has a Lipschitz nonlinearity, and its viscosity solution
is given by a backward stochastic differential equation
\cite{PardouxPeng1992,EKPQ1997} (Proposition~\ref{prop:rp-bsde-identification}).
The field $Y$ is then obtained by linearizing the associated forward--backward
system with respect to the intensity. The differentiability of solutions of
backward equations with respect to a parameter is well known (see \cite{EKPQ1997}
and, for generators of quadratic growth, \cite{AIR2007}). In our setting, however,
the parameter enters through the drift of the component $y$, which has no diffusion,
and we need an estimate that is uniform in the state and in time. In
Section~\ref{sec:rp-numerics} we use a numerical example to illustrate the orders
in \eqref{eq:intro-main} and in Corollary~\ref{theorem:rp-mean-variance}.

The paper is organized as follows. Section~\ref{sec:results} treats the
$p$-Laplacian equation: Section~\ref{subsec:plap} contains the formulation and
preliminaries, and Section~\ref{subsec:reg} the regularized problem and its
estimates. In Section~\ref{sec:model} we recall the two-factor model and
introduce the random perturbation of the prepayment intensity. The expansion
\eqref{eq:intro-main} and its consequences are stated in
Section~\ref{subsec:rp}, and the numerical example is in
Section~\ref{sec:rp-numerics}. The proofs of the results of
Section~\ref{subsec:rp} are given in Section~\ref{sec:proofs}.

\section{Deterministic case with p-Laplacian diffusion}\label{sec:results}

\subsection{Analysis of the p-Laplacian case}\label{subsec:plap}
We consider the p-laplacian analog of the MBS equation. There are several motivations for this from the literature, \cite{Takac2023, parshad2026, li2020}.

\begin{equation}\label{eq:rp-pricing-pde-plap}
\begin{cases}
\displaystyle u_t=a \nabla \cdot \left(|\nabla u|^{p-2}\nabla u\right)+c \cdot \nabla u
 -k\dfrac{|\nabla u|^2}{u+h+\xi}+(\tau-r)h-ru,
 & (z,t)\in \Omega \times(0,T),\\[1ex]
 u(z,0)=g(z), z\in  \Omega \subset \mathbb{R}^{n}, n=1,2, \ |\Omega| < \infty,  \ 1<p<2. \\
  \nabla u \cdot \textbf{n} = 0 , z \in \partial \Omega.
\end{cases}
\end{equation}

\subsubsection{Functional Preliminaries}

\begin{lemma}
\label{lem:wlp}
Consider exponents $0 < p_0 < p_1 < \infty$ and a domain $\Omega$ that is a closed and bounded in $\mathbb{R}^{n}$, $n\geq 1$, and a $u(x) \in L^{p_1}(\Omega)$. Then,

\begin{equation}
||u||_{L^{p_{\theta}}(\Omega)} \leq C||u||^{1-\theta}_{L^{p_{0}}(\Omega)} ||u||^{\theta}_{L^{p_{1}}(\Omega)}
\end{equation}

for all $0 \leq \theta \leq 1$, and where $p_{\theta}$ is defined as $\frac{1}{p_{\theta}} = \frac{1-\theta}{p_0} + \frac{\theta}{p_1}$.
\end{lemma}

\begin{lemma}
\label{lem:gns}
Consider a $\phi \in W^{m,q^{'}}(\Omega) \cap L^{q}(\Omega)$. Then if $p^{'}, q^{'}, q \geq 1, 0 \leq \theta \leq 1$, and

\begin{equation}
k - \frac{n}{p^{'}} \leq \theta \left(  m - \frac{n}{q^{'}} \right) - (1 - \theta) \frac{n}{q},
\end{equation}
there exists a constant C such that,

\begin{equation}
	||\phi||_{W^{k,p^{'}}(\Omega)} \leq C ||\phi||^{\theta}_{W^{m,q^{'}}(\Omega)} ||\phi||^{1 - \theta}_{L^{q}(\Omega)} 
	\end{equation}

\end{lemma}	

Now we state the classical Aubin-Lions compactness Lemma,

\begin{lemma}
\label{lem:al}
Let $X_0, X$ and $X_1$ be three reflexive Banach spaces with 

$X_0 \hookrightarrow \hookrightarrow X \hookrightarrow X_1$. That is $ X_0$ is compactly embedded in $X$ and that $X$ is continuously embedded in $X_1$. For $1\leq p,q \leq \infty$, let

 \begin{equation}
W = \{ u \in L^{p^{'}}([0,T];X_0) \ | \  u^{'} \in L^{q^{'}}([0,T];X_1)  \}
\end{equation}

(i) If $p^{'} < \infty$ then the embedding of $W$ into $L^{p^{'}}([0,T];X)$ is compact.

(ii) If $p^{'} = \infty$ and $q^{'} > 1$ then the embedding of $W$ into $C([0,T];X)$ is compact.
\end{lemma}

\begin{lemma}
\label{lem:pc1}
    Consider functions $\phi, \sigma \in C^{1}(\Omega), \Omega \subset \mathbb{R}^{n}, n=1,2, |\Omega| < \infty$, $\epsilon > 0$. Then if $p \geq 2$ we have,
    \begin{equation}
        (|\nabla \phi|^{2} + \epsilon)^{\frac{p-2}{2}}|\nabla \phi|^{2} \geq |\nabla \phi|^{p},
\end{equation}

\begin{equation}
        (|\nabla \phi|^{2} + \epsilon)^{\frac{p-2}{2}}|\nabla \phi| \leq C\left(|\nabla \phi|^{p}+1\right).
\end{equation}

and 

\begin{equation}
        (|\phi|^{p-2}\phi - |\sigma|^{p-2}\sigma)\cdot (\phi - \sigma) \geq C |\phi - \sigma|^{p}
\end{equation}

\end{lemma}

\begin{lemma}
\label{lem:pc1aa}
    Consider functions $\phi, \sigma \in C^{1}(\Omega), \Omega \subset \mathbb{R}^{n}, n=1,2, |\Omega| < \infty$, $\epsilon > 0$. Then if $1< p <2$ we have,
    \begin{equation}
        (|\nabla \phi|^{2} + \epsilon)^{\frac{p-2}{2}}|\nabla \phi|^{2} \geq \left(|\nabla \phi|^{2}+\epsilon \right)^{\frac{p}{2}} - (\epsilon)^{\frac{p}{2}} \,
\end{equation}

\end{lemma}

\subsubsection{Weak Solutions}

\begin{definition}
\label{def:d11}
Let $1<p<2$, let $T \in (0,\infty]$ and $\Omega \subset \mathbb{R}^{n}, n=1,2$ be a bounded domain, with smooth boundary. A non-negative function $u$ defined on $\Omega \times (0,T)$ is called a weak solution to \eqref{eq:rp-pricing-pde-plap} on $[0,T]$ if,

(i) $u  \in L^{2}_{loc}([0,T];L^{2}(\Omega))$

(ii) $|\nabla u|^{p-2} \nabla u \in L^{1}_{loc}([0,T];L^{1}(\Omega));$

(iii) For every $\phi \in C^{\infty}_{0}(\Omega \times [0,T))$ we have,

\begin{eqnarray}
&& \int^{T}_{0}\int_{\Omega} u \phi_{t} dxdt - \int_{\Omega} u_{0}(x) \phi(x,0)dx \nonumber \\
&& = 
- \int^{T}_{0}\int_{\Omega} |\nabla u|^{p-2} \nabla u \cdot \nabla \phi dx dt + \int^{T}_{0}\int_{\Omega}\textbf{c} u \nabla \phi dx dt + \int^{T}_{0}\int_{\Omega} \left( k\dfrac{|\nabla u|^2}{u+h+\xi}+(\tau-r)h-ru\right) \phi dxdt . \nonumber \\
\end{eqnarray}
\end{definition}

\subsection{The regularized system and preliminary estimates}\label{subsec:reg}
There can be several technical challenges with \eqref{eq:rp-pricing-pde-plap}, and well-posedness has to be dealt with. We follow many of the techniques enunciated from \cite{li2020, parshad2026}, in the ensuing analysis. Also note, in the ensuing analysis $C,C_{i}$ are generic constants and can change in value from line to line. We first posit the regularized system, in a bounded spatial domain with dimension $n \leq 2$.

\begin{equation}\label{eq:rp-pricing-pde-plap-reg}
\begin{cases}
\displaystyle u^{\epsilon}_t=a \nabla \cdot \left(|\nabla u^{\epsilon}|^{p-2}+\epsilon)\nabla u^{\epsilon}\right)+ \nabla  \cdot (\textbf{c} u^{\epsilon})
 -k\dfrac{|\nabla u^{\epsilon}|^2}{u^{\epsilon}+h+\xi}+(\tau-r)h-ru^{\epsilon},
 & (z,t)\in \Omega \times(0,T),\\[1ex]
\nabla u^{\epsilon} \cdot \eta=0, z \in \partial \Omega, \  u^{\epsilon}(z,0)=g(z), z\in  \Omega, \ 1<p<2.
\end{cases}
\end{equation}

We state the following theorem, about the regularized system,

\begin{theorem}
\label{thm:mt11}
    Consider \eqref{eq:rp-pricing-pde-plap-reg}. Let $p>\frac{3}{2}$. If the initial data $u_{0} \in L^{\infty}(\Omega)$, is sufficiently small, that is $||u_{0}||_{\infty} < C_{1}$, then for any  $q \geq 1$, there exists a $C=C(q)$ such that for any $ \epsilon > 0$, we have,
\begin{equation}
    ||u^{\epsilon}||_{L^{q}(\Omega)} \leq C,
\end{equation}
  for all $t \in [0, T_{max, \epsilon})$.
\end{theorem}

\begin{proof}
We multiply the equation for $u^{\epsilon}$ in \eqref{eq:rp-pricing-pde-plap-reg} by $q(u^{\epsilon})^{q-1}$, to obtain,

\begin{eqnarray}
&& \frac{d}{dt}||u^{\epsilon}||^{q}_{q} +   q(q-1)\int_{\Omega} \left(
|\nabla u^{\epsilon}|^{2}+{\epsilon}\right)^{\frac{p-2}{2}}|\nabla u^{\epsilon}|^{2} (u^{\epsilon})^{q-2} dx \nonumber \\
&& -\int_{\partial \Omega} \left(
|\nabla u^{\epsilon}|^{2}+{\epsilon}\right)^{\frac{p-2}{2}}( \nabla u^{\epsilon} \cdot \eta ) ((qu^{\epsilon})^{q-1}) dS \nonumber \\
&&+r||u^{\epsilon}||^{q}_{q} 
 + k\int_{\Omega} \dfrac{|\nabla u^{\epsilon}|^2 (u^{\epsilon})^{q-1}}{u^{\epsilon}+h+\xi}dx \nonumber \\
 && = q(q-1)\int_{\Omega} (\textbf{c}\cdot \nabla u^{\epsilon}) (u^{\epsilon})^{q-1}dx + \int_{\Omega}(\tau-r)h (u^{\epsilon})^{q-1} dx \nonumber \\
 \nonumber \\
\end{eqnarray}

Note via the boundary conditions prescribed we have,

\begin{eqnarray}
&& \frac{d}{dt}||u^{\epsilon}||^{q}_{q} +   q(q-1)\int_{\Omega} \left(
|\nabla u^{\epsilon}|^{2}+{\epsilon}\right)^{\frac{p-2}{2}}|\nabla u^{\epsilon}|^{2} (u^{\epsilon})^{q-2} dx \nonumber \\
&&+r||u^{\epsilon}||^{q}_{q} 
 + k\int_{\Omega} \dfrac{|\nabla u^{\epsilon}|^2 (u^{\epsilon})^{q-1}}{u^{\epsilon}+h+\xi}dx \nonumber \\
 && = q(q-1)\int_{\Omega} (\textbf{c}\cdot \nabla u^{\epsilon}) (u^{\epsilon})^{q-1}dx + \int_{\Omega}(\tau-r)h (u^{\epsilon})^{q-1} dx \nonumber \\
 \nonumber \\
\end{eqnarray}

Now note when $1<p<2$, we have,
\begin{equation}
\label{eq:in1}
\int_{\Omega} \left(
|\nabla u^{\epsilon}|^{2}+{\epsilon}\right)^{\frac{p-2}{2}}|\nabla u^{\epsilon}|^{2} (u^{\epsilon})^{q-2} dx
\geq 
\int_{\Omega} \left(
|\nabla u^{\epsilon}|\right)^{p}(u^{\epsilon})^{q-2} dx - \epsilon^{\frac{p}{2}}\int_{\Omega} (u^{\epsilon})^{q-2} dx .
\end{equation}
Also note,
\begin{equation}
\label{eq:in2}
\int_{\Omega} (u^{\epsilon})^{q-2}|\nabla u^{\epsilon}|^{p}dx = \frac{p^{p}}{(p+q-2)^{p}}||\nabla (u^{\epsilon})^{\frac{p+q-2}{p}}||^{p}_{p}.
\end{equation}
Using the above in the $\int_{\Omega} (\textbf{c}\cdot \nabla v^{\epsilon}) (v^{\epsilon})^{q-1}dx$ term, by setting $p=1$, and $q=q^{'}+1$, where $q^{'}$, is a dummy variable, yields, 

\begin{equation}
\int_{\Omega} (u^{\epsilon})^{q-1}|\nabla u^{\epsilon}|^{1}dx = \frac{1}{q}||\nabla (u^{\epsilon})^{q}||^{1}_{1}.
\end{equation}
Using positivity and the earlier estimates via \eqref{eq:in1}-\eqref{eq:in2} entails,

\begin{eqnarray}
&& \frac{d}{dt}||u^{\epsilon}||^{q}_{q}  +  q(q-1)[\int_{\Omega} \left(
|\nabla u^{\epsilon}|\right)^{p}(u^{\epsilon})^{q-2} dx - \epsilon^{\frac{p}{2}}\int_{\Omega} (u^{\epsilon})^{q-2} dx]
 \nonumber \\
&& 
+r||u^{\epsilon}||^{q}_{q} 
 + k\int_{\Omega} \dfrac{|\nabla u^{\epsilon}|^2 (u^{\epsilon})^{q-1}}{u^{\epsilon}+h+\xi}dx  
 \nonumber \\
&& 
\leq \frac{1}{q}||\nabla (u^{\epsilon})^{q}||^{1}_{1} + \int_{\Omega}(\tau-r)h (u^{\epsilon})^{q-1} dx\nonumber \\
\end{eqnarray}
Again, using positivity we have,

\begin{eqnarray}
&& \frac{d}{dt}||u^{\epsilon}||^{q}_{q}  +  q(q-1)\int_{\Omega} \left(
|\nabla u^{\epsilon}|\right)^{p}(u^{\epsilon})^{q-2} dx
 \nonumber \\
&& 
+r||u^{\epsilon}||^{q}_{q} 
 + k\int_{\Omega} \dfrac{|\nabla u^{\epsilon}|^2 (u^{\epsilon})^{q-1}}{u^{\epsilon}+h+\xi}dx  
 \nonumber \\
&& 
\leq \frac{1}{q}||\nabla (u^{\epsilon})^{q}||^{1}_{1} + \int_{\Omega}(\tau-r)h (u^{\epsilon})^{q-1} dx + \epsilon^{\frac{p}{2}}\int_{\Omega} (u^{\epsilon})^{q-2} dx\nonumber \\
\end{eqnarray}

Notice,
\begin{equation}
q(q-1)\int_{\Omega} (\textbf{c}\cdot \nabla u^{\epsilon}) (u^{\epsilon})^{q-1}dx = q(q-1)\int_{\Omega} (\textbf{c}\cdot \nabla u^{\epsilon}) (u^{\epsilon})^{\frac{q-2}{p}} (u^{\epsilon})^{q^{*}} ,
\end{equation}
where 
\begin{equation}
    q^{*} = \left(q-1 - \frac{q-2}{p}\right)= \frac{p(q-1) -(q-2)}{p}.
\end{equation}
Now, using Holder-Young with $\epsilon_{1}$, and exponents $p, q=\frac{p}{p-1}$, we have,

\begin{equation}
\int_{\Omega} (u^{\epsilon})^{q-1}|\nabla u^{\epsilon}|^{1}dx
\leq \frac{\epsilon_{1}}{p}  q(q-1)\int_{\Omega} (u^{\epsilon})^{q-2}|\nabla u^{\epsilon}|^{p}dx + \frac{p-1}{p}f(\epsilon_{1})\int_{\Omega} (u^{\epsilon})^{q^{**}}dx.
\end{equation}
Here, $q^{**} = q^{*}\left(\frac{p}{p-1}\right)$, so,

\begin{eqnarray}
   && q^{**} \nonumber \\
   && = \left(q-1 - \frac{q-2}{p}\right)\left( \frac{p}{p-1}\right)  \nonumber \\
   && = \frac{p(q-1) -(q-2)}{p-1} = \frac{(p-1)q + (2-p)}{p-1}  \nonumber \\
   && = q + \left(\frac{2-p}{p-1}\right).  \nonumber \\
\end{eqnarray}
Choosing $\epsilon_{1}$, such that $\epsilon_{1} < \frac{p}{2} $, and inserting this into earlier estimates we have,

\begin{eqnarray}
\label{eq:in44}
&& \frac{d}{dt}||u^{\epsilon}||^{q}_{q} + r||u^{\epsilon}||^{q}_{q} + \frac{q}{2}(q-1)\int_{\Omega} \left(|\nabla u^{\epsilon}|^{p}\right)(u^{\epsilon})^{q-2} dx 
 + k\int_{\Omega} \dfrac{|\nabla u^{\epsilon}|^2 (u^{\epsilon})^{q-1}}{u^{\epsilon}+h+\xi}dx  
 \nonumber \\
&& \leq  f(\epsilon_{1})\int_{\Omega} (u^{\epsilon})^{q^{**}}dx + \int_{\Omega}(\tau-r)h (u^{\epsilon})^{q-1} dx + \epsilon^{\frac{p}{2}}\int_{\Omega} (u^{\epsilon})^{q-2} dx. \nonumber \\
\end{eqnarray}
Now, note that 
\begin{equation}
    q^{**} = q + \left(\frac{2-p}{p-1}\right) < q+1, 
\end{equation} 
 if $p > \frac{3}{2}$,  then using the embedding of $L^{q+1}(\Omega) \hookrightarrow L^{q^{**}}(\Omega)$ entails,

\begin{eqnarray}
\label{eq:in4411}
&& \frac{d}{dt}||u^{\epsilon}||^{q}_{q} + r||u^{\epsilon}||^{q}_{q} + \frac{q}{2}(q-1)\int_{\Omega} \left(|\nabla u^{\epsilon}|^{p}\right)(u^{\epsilon})^{q-2} dx 
 + k\int_{\Omega} \dfrac{|\nabla u^{\epsilon}|^2 (u^{\epsilon})^{q-1}}{u^{\epsilon}+h+\xi}dx  
 \nonumber \\
&& \leq  f(\epsilon_{1})\int_{\Omega} (u^{\epsilon})^{q+1}dx + \int_{\Omega}(\tau-r)h (u^{\epsilon})^{q-1} dx + \epsilon^{\frac{p}{2}}\int_{\Omega} (u^{\epsilon})^{q-2} dx. \nonumber \\
\end{eqnarray}

Now we use the embedding of $L^{q}(\Omega) \hookrightarrow L^{q-1}(\Omega) \hookrightarrow L^{q-2}(\Omega)$  entails comparison to the ODE, $y^{'}=C_{1} - C_{2}y + C_{3}y^{1+\frac{1}{q}}$, for $C_{1}, C_{2}, C_{3}, q > 0$.
This enables a bound on $y$, for small enough initial data, via standard ODE comparison.

\end{proof}
\begin{remark}
    Note,  $\boxed{q^{**}=q + \left(\frac{2-p}{p-1}\right)}$, as $p \searrow 1$, $q^{**} >>1$, is arbitrarily large> However, for $p>\frac{3}{2}$, this reduces to $q^{**} < q+1$, and hence
    $||v^{\epsilon}||^{q+1}_{q+1}$, which enables an ODE comparison for sufficiently small initial data.
\end{remark}

\begin{lemma}
\label{lem:ge11n}

 Consider \eqref{eq:rp-pricing-pde-plap-reg}. Let $p>\frac{3}{2}$. If the initial data $u_{0} \in L^{\infty}(\Omega)$, is sufficiently small, that is $||u_{0}||_{\infty} < C_{1}$, then for any  $q \geq 1$, there exists a $C$ depending only on the problem parameters, such that we have,

\begin{equation}
\label{eq:12l}
    ||u^{\epsilon}||_{L^{\infty}(\Omega)} \leq C,
\end{equation}

\begin{equation}
\label{eq:13l}
\int^{t}_{0}\int_{\Omega} 
|\nabla u^{\epsilon}(.,s)|^{p}(u^{\epsilon}(.,s))^{q-2} dx ds \leq C
\end{equation}

\begin{equation}
\label{eq:14l}
\int^{t}_{0}\int_{\Omega} 
|\nabla u^{\epsilon}(.,s)|^{p} dx ds \leq C
\end{equation}
\end{lemma}

\begin{proof}
   The proof follows via similar estimates in \cite{li2020, parshad2026}, as well as 
    integration in time $[0,T]$ of the inequalities in \eqref{eq:in4411}, followed by Moser-Alikakos type iteration, as in \cite{li2020}. This proves the lemma.
\end{proof}

\begin{theorem}
\label{thm:vubd}
     Consider \eqref{eq:rp-pricing-pde-plap-reg}. Let $p>\frac{3}{2}$. If the initial data $u_{0} \in L^{\infty}(\Omega)$, is sufficiently small, that is $||u_{0}||_{\infty} < C_{1}$, then 
 there exists a function $u \in L^{p}_{loc}((0,\infty), W^{1,p}(\Omega)) \cap L^{\infty}((0,\infty);L^{\infty}(\Omega))$   as well as a $\Gamma \in L^{\frac{p}{p-1}}_{loc}((0,\infty);L^{\frac{p}{p-1}}(\Omega)) \cap L^{\infty}((0,\infty);L^{\infty}(\Omega))$, and a sequence of approximants $\epsilon = \epsilon_{j} \searrow 0$, such that,

 \begin{equation}
u^{\epsilon} \overset{*}{\rightharpoonup} u^{*}  \ \mbox{in} \ L^{\infty}((0,\infty);L^{\infty}(\Omega))
\end{equation}

 \begin{equation}
 \label{eq:w1pv}
\nabla u^{\epsilon} \rightharpoonup \nabla u^{*}  \ \mbox{in} \ L^{p}_{loc}((0,\infty);L^{p}(\Omega))
\end{equation}

 \begin{equation}
|\nabla u^{\epsilon}|^{p-2}\nabla u^{\epsilon} \rightharpoonup  \Gamma \ \mbox{in} \ L^{\frac{p}{p-1}}_{loc}((0,\infty);L^{\frac{p}{p-1}}(\Omega))
\end{equation}

\end{theorem}

\begin{proof}
    The proof follows via the uniform bounds in Lemma \ref{lem:ge11n}.
\end{proof}

Next, we state a result about the time derivative of the solution,

\begin{theorem}
\label{thm:w1s}
    Consider \eqref{eq:rp-pricing-pde-plap-reg}. Let $2>p>\frac{3}{2}$. If the initial data $u_{0}(x) \in L^{\infty}(\Omega)$, then there exists a constant $C$ such that for any time $T > 0$, we have,

     \begin{equation} 
\left \lVert \frac{\partial u^{\epsilon}}{\partial t}\right \rVert_{L^{1}((0,T); (W^{1,p}(\Omega))^{*})} \leq C
\end{equation}

\end{theorem}

\begin{proof}
The proof follows via standard methods, see \cite{li2020, parshad2026}.
\end{proof}

\begin{lemma}
\label{lem:vlp}
Consider \eqref{eq:rp-pricing-pde-plap-reg}, with $ 1 < p \leq 2$, there exists a subsequence ${\epsilon_{j}}$, such that,
 \begin{equation}
u^{\epsilon_{j}} \rightarrow u^{*} \ \mbox{in} \  L^{p}([0,T];L^{p}(\Omega))
\end{equation}

\end{lemma}

\begin{proof}
    We know via Lemma \ref{lem:ge11n} that we have,
    \begin{equation}
u^{\epsilon_{j}} \in L^{p}((0,T), W^{1,p}(\Omega)).
\end{equation}

Furthermore, via Theorem \ref{thm:w1s} we have,

 \begin{equation}
 \frac{\partial u^{\epsilon_{j}}}{\partial t} \in L^{1}((0,T); (W^{1,p}(\Omega))^{*}).
\end{equation}
Thus the classical Aubin-Lion compactness Lemma yields,

 \begin{equation}
W^{1,p}(\Omega) \hookrightarrow \hookrightarrow L^{p}(\Omega) \hookrightarrow \left(W^{1,p}(\Omega)\right)^{*},
 \end{equation}
 to obtain $u^{\epsilon_{j}} \rightarrow u^{*}$ in 
 $L^{p}((0,T), L^{p}(\Omega))$.

We derive, 
\begin{equation}
\label{eq:56e}
\lim_{j \rightarrow \infty}  \int^{T}_{0}\int_{\Omega}|u^{\epsilon_{j}}-u^{*}|dxdt \rightarrow 0 .
 \end{equation}
Note, as $\lim_{j \rightarrow \infty}\epsilon_{j} \rightarrow 0$. This follows via 
\begin{equation}
L^{p}((0,T), L^{p}(\Omega)) \hookrightarrow L^{1}((0,T), L^{1}(\Omega)).
 \end{equation}
this proves the lemma.
\end{proof}

\begin{lemma}
\label{lem:cnt}
Consider \eqref{eq:rp-pricing-pde-plap-reg} with $1<p<2$. There exists a subsequence ${\epsilon_{j}}$, such that for any $T>0$,

\begin{equation}
|\nabla u^{\epsilon_{j}}|^{p-2}  \nabla u^{\epsilon_{j}} \rightarrow |\nabla u|^{p-2}  \nabla u \ \mbox{in} \ L^{\frac{p}{p-1}}(\Omega \times (0,T)).
 \end{equation}
\end{lemma}
\begin{proof}
    The proof follows via similar estimates in \cite{li2020, parshad2026}.
\end{proof}

We can now state our main result,

\begin{theorem}
\label{thm:mt1}
Consider the system \eqref{eq:rp-pricing-pde-plap}, with $2>p>\frac{3}{2}, >0$. If the initial data $u_{0}(x) \in L^{\infty}(\Omega)$, then there exists a global weak solution $u$ to system \eqref{eq:rp-pricing-pde-plap}, in the sense of Definition \ref{def:d11}. Furthermore, there exists a constant $C$, independent of time and initial data, such that

\begin{equation}
   ||u||_{L^{\infty}(\Omega)} \leq C, \
  \end{equation}
  for all $t > 0$.
\end{theorem}

\begin{proof}
The proof follows via Lemma \ref{lem:ge11n}, Theorem \ref{thm:vubd}, Theorem \ref{thm:w1s}, Lemma \ref{lem:vlp} and Lemma \ref{lem:cnt}.
\end{proof}

\section{Stochastic case with random environment}\label{sec:model}
Let $t=T-s$ be time to maturity and $z=(x,y)$, where $x$ is a financial factor and $y$ is cumulative prepayment. As in \cite[Example~10.26]{PapiBriani2011},the state has dynamic 
\begin{equation}\label{eq:rp-state-dynamic}
    dX_s=b(X_s,T-s)\,ds+\sigma(T-s)\,dB_s\qquad
    d\mathcal A_s=\eta(X_s,\mathcal A_s, T-s)\,ds
\end{equation}

Let $\mathit{MB}(s)$ denote scheduled remaining principal without prepayment. For example, a fixed mortgage rate $\tau'>0$ and an fixed initial balance $\mathit{MB}(0)>0$ give \[
\mathit{MB}(s)=\mathit{MB}(0)
\frac{e^{\tau' T}-e^{\tau' s}}{e^{\tau' T}-1}
\]
Let $\tau>0$ be the coupon rate, $r\ge0$ the constant discount rate, $A_o>0$ the initial value of the risk-free asset, and $\rho>0$ the risk parameter. Define 
\begin{equation}\label{eq:rp-balance-and-price}
 h(y,t)=\mathit{MB}(T-t)e^{-s_0(y)},\qquad
 \xi(t)=A_0e^{r(T-t)},\qquad P=u+h,
\end{equation}
where $s_0(y)=y$ for $y\ge0$, extended smoothly to $y<0$ so that the pricing equation can be
posed on $Q_T=\mathbb R^2\times[0,T]$. Thus $h$ is the remaining principal, $P$ is the full security value and $u$ is its adjustment. With $a(t)=\sigma(t)^2/2$ and $k(t)=\rho\sigma(t)^2$, the two state specialization of the pricing equation in \cite[(10.50)]{PapiBriani2011} is
\begin{equation}\label{eq:rp-pricing-pde}
\begin{cases}
\displaystyle u_t=a u_{xx}+b u_x+\eta u_y
 -k\dfrac{u_x^2}{u+h+\xi}+(\tau-r)h-ru,
 & (z,t)\in\mathbb R^2\times(0,T),\\[1ex]
 u(z,0)=0,&z\in\mathbb R^2.
\end{cases}
\end{equation}

Let $\eta^0\ge0$ be the baseline intensity, and put
\begin{equation}\label{eq:re-scenario-set}
 \Theta=\{\theta\in C([0,T]):\|\theta\|_\infty\le\tfrac12\},
 \qquad\eta^\theta(z,t)=\eta^0(z,t)(1+\theta(t)).
\end{equation}

Write $u^\theta$ for the solution with intensity $\eta^\theta$(see Proposition~\ref{prop:rp-bsde-identification}) and $u^0=u^{\theta=0}$. On $(\Omega,\mathcal{G},\mathbb P)$
Let $\zeta$ be a $C([0,T])$-valued random variable and $K>0$ a constant such that
\begin{equation}\label{eq:rp-random-intensity}
N(\omega)=\|\zeta(\cdot,\omega)|_\infty\le K\quad\text{a.s.},\quad\eta^\varepsilon(z,t,\omega)=\eta^0(z,t)(1+\varepsilon\zeta(t,\omega))
\end{equation}

For $|\varepsilon|\le(2K)^{-1}$ we have $\|\varepsilon\zeta(\cdot,\omega)\|_\infty\le \tfrac12$,
so $\varepsilon\zeta(\cdot,\omega)\in\Theta$ for a.e\ $\omega$, and we may put
$u^\varepsilon(\cdot,\omega)=u^{\varepsilon\zeta(\omega)}$ and $P^\varepsilon=u^\varepsilon+h$

\section{Random prepayment at $p=2$}\label{subsec:rp}
In this subsection we return to the Papi-Briani equation \eqref{eq:rp-pricing-pde}, which has diffusion in $x$ only and is posed on $Q_T$, and perturb the prepayment intensity by the random factor $1+\varepsilon\zeta$ \eqref{eq:rp-random-intensity}. Our main result, Theorem~\ref{theorem:rp-first-order}, shows that $u^\varepsilon-u^0$ equals $\varepsilon Y$ up to an error of order $\varepsilon^2$, uniformly on $Q_T$, where $Y$ is linear in $\zeta$. The mean, the variance and the mean-square error of $P^\varepsilon$ then follow by taking expectations (Corollary~\ref{theorem:rp-mean-variance}). The proof uses a backward SDE representation of $u^\theta$, set up in Proposition~\ref{prop:rp-bsde-identification}. All proofs are given in Section~\ref{sec:proofs}.
We make the following assumptions
\begin{enumerate}
\item[(A1)] \emph{Deterministic data.}
$T>0$, $A_0>0$, $\rho>0$, and $\tau\ge r\ge0$ with $\tau>0$;
$\mathit{MB}\in C^1([0,T])$, $\mathit{MB}\ge0$;
$s_0\in C^3(\mathbb R)$, $s_0(y)=y$ for $y\ge0$,
$\inf s_0> -\infty$, and $s_0',s_0'',s_0'''$ are bounded;
$\sigma\in C([0,T])$ is nonnegative;
$b\in C(\mathbb R\times[0,T])$ is bounded and globally Lipschitz in
$x$, uniformly in $t$;
$\eta^0\in C(\mathbb R^2\times[0,T])$ is nonnegative and
$\sup_{t\in[0,T]}\|\eta^0(\cdot,t)\|_{W^{2,\infty}(\mathbb R^2)}<\infty$.
\item[(A2)] \emph{Random environment.} $\zeta$ is as in
\eqref{eq:rp-random-intensity}.
\end{enumerate}

Write $\|\cdot\|_\infty$ for the supremum norm on $Q_T$.  Let $\nu$ be
any Borel probability measure on $\mathbb R^2$, and set
\begin{equation}\label{eq:rp-spaces}
 H=L^2(\mathbb R^2,\nu),\qquad
 \mathcal X=C([0,T];H),\qquad
 \|v\|_{\mathcal X}=\sup_{0\le t\le T}\|v(t)\|_H.
\end{equation}
For bounded continuous $v$, $\|v\|_{\mathcal X}\le\|v\|_\infty$.

\begin{theorem}[Uniform first-order sensitivity]
\label{theorem:rp-first-order}
Assume (A1).  There exists a bounded deterministic linear operator
$\mathcal L:C([0,T])\to C_b(Q_T)$ and a constant $C_2<\infty$,
depending only on (A1), such that for all $\theta\in\Theta$,
\begin{equation}\label{eq:rp-Frechet}
 \|u^\theta-u^0-\mathcal L\theta\|_\infty
 \le C_2\|\theta\|_\infty^2.
\end{equation}
In particular, under (A2), $Y=\mathcal L\zeta$ is strongly measurable
as an $\mathcal X$-valued random variable.  For
$|\varepsilon|\le\varepsilon_*:=(2K)^{-1}$, almost surely,
\begin{equation}\label{eq:rp-main}
 \|u^\varepsilon-u^0-\varepsilon Y\|_\infty
 \le C_2\varepsilon^2N(\omega)^2.
\end{equation}
Writing $R^\varepsilon=u^\varepsilon-u^0-\varepsilon Y$, one has
\begin{equation}\label{eq:rp-main-L2}
 \|R^\varepsilon\|_{L^2(\Omega;\mathcal X)}
 \le C_2\varepsilon^2\big(\mathbb EN^4\big)^{1/2}.
\end{equation}
Consequently, $\theta\mapsto u^\theta$ is Fr\'echet differentiable
at zero in the supremum norm, and $\varepsilon\mapsto u^\varepsilon$
is differentiable at zero in $L^2(\Omega;\mathcal X)$ with derivative
$Y$.  The same expansion holds for $P^\varepsilon$ since $h$ does not
depend on $\varepsilon$.
\end{theorem}

The operator $\mathcal{L}$ is deterministic, so the first-order response depends on the randomness only through $\zeta$, via a fixed linear map. In Section~\ref{sec:proofs}, $\mathcal{L}$ is given explicitly by a linearized forward-backward system.

Since the remainder in \eqref{eq:rp-main} is bounded by $C_2K^2\varepsilon^2$ almost surely, we may take expectations in the expansion. This gives the following consequences for the mean and the variance of the price.

\begin{corollary}[Mean-square sensitivity, mean bias and variance]
\label{theorem:rp-mean-variance}
Assume (A1)--(A2).  There exists $C_3<\infty$ such that for all
$|\varepsilon|\le\varepsilon_*$ and $t\in[0,T]$,
\begin{align}
 \Big|\mathbb E\|P^\varepsilon(t)-P^0(t)\|_H^2
       -\varepsilon^2\mathbb E\|Y(t)\|_H^2\Big|
 &\le C_3|\varepsilon|^3,\label{eq:rp-ms}\\
 \|\mathbb EP^\varepsilon(t)-P^0(t)-\varepsilon\mathbb EY(t)\|_H
 &\le C_2K^2\varepsilon^2,\label{eq:rp-mean}\\
 \big|\operatorname{Var}_H(P^\varepsilon(t))
       -\varepsilon^2\operatorname{Var}_H(Y(t))\big|
 &\le C_3|\varepsilon|^3,\label{eq:rp-var}
\end{align}
where $\operatorname{Var}_H(V)=\mathbb E\|V-\mathbb EV\|_H^2$.

\end{corollary}

Thus the mean price moves at first order by $\varepsilon\,\mathcal{L}(\mathbb E \zeta)$, which vanishes for centered noise, and the variance is $\varepsilon^2\operatorname{Var}_H(Y(t))$ to leading order. Since $Y$ is linear in $\zeta$, this leading coefficient depends on the law of $\zeta$ only through its covariance function.

We now set up the representation used in the proofs. For $\theta\in\Theta$ write $D^\theta=u^\theta+h+\xi$, and set 
\begin{equation}\label{eq:rp-transform}
    \Phi(d)=\int_{A_0}^d q^{-2\rho}\,dq,\qquad v^\theta=\Phi(D^\theta).
\end{equation}
Since $k=2 \rho a$, we have $-\tfrac{k}{d}\Phi'(d)-a\Phi''(d)=0$, so this change of unknown removes the quadratic gradient term in \eqref{eq:rp-pricing-pde}. Put 
\begin{equation}\label{eq:rp-price-barrier}
    K_0=\sup_{(z,t)\in Q_T}\frac{(\tau-r)h(z,t)}{1+rt} \qquad
    M_D=K_0T+\|h\|_\infty+\|\xi\|_\infty
    \end{equation}
and $I=\Phi([A_0,M_D])$. By Proposition~\ref{prop:rp-bsde-identification} below,
$A_0\le D^\theta\le M_D$, so $v^\theta$ takes values in $I$.Fix a smooth $\psi:\mathbb R\to(0,\infty)$ that equals $\Phi^{-1}$ on a neighborhood of $I$ and is constant outside a compact interval, and let 
\begin{align}
G(y)&=\Phi\bigl(h(y,0)+\xi(0)\bigr),\label{eq:rp-G}\\
f^\theta(z,t,w)&=\psi(w)^{-2\rho}\bigl[\tau h(y,t)+h_t(y,t)-\eta^\theta(z,t)h_y(y,t)-r\psi(w)\bigr].\label{eq:rp-transformed-driver}
\end{align}

Using $h_x=0$ and $\xi'=-r\xi$, one checks that $v^\theta$ solves
$v_t=av_{xx}+bv_x+\eta^\theta v_y+f^\theta(z,t,v)$ with $v(\cdot,0)=G$. Since $\psi=\Phi^{-1}$
on the range of $v^\theta$, the cutoff will not change this equation, and it makes $f^\theta$
bounded and globally Lipschitz in $(z,w)$.

Let $(\widehat\Omega,\widehat{\mathcal F},\widehat{\mathbb P})$ carry a one-dimensional Browinian motion $B$ with its completed, right continuous natural filtration, and write
$\widehat{\mathbb E}$ for the expectation. For $(z,t)\in Q_T$ and $\theta\in\Theta$, let
$S^\theta_s=(X_s,\mathcal A^\theta_s)$, $0\le s\le t$, start at $z$ and solve
\begin{equation}\label{eq:rp-forward}
 dX_s=b(X_s,t-s)\,ds+\sigma(t-s)\,dB_s,\qquad
 d\mathcal A^\theta_s=\eta^\theta(S^\theta_s,t-s)\,ds,
\end{equation}
and let $(V^\theta,Q^\theta)$ solve
\begin{equation}\label{eq:rp-bsde}
 V^\theta_s=G(\mathcal A^\theta_t)+\int_s^t f^\theta(S^\theta_r,t-r,V^\theta_r)\,dr
 -\int_s^t Q^\theta_r\,dB_r.
\end{equation}

We write $V^{\theta,z,t}$ when the starting point matters.

\begin{proposition}
\label{prop:rp-bsde-identification}
Assume (A1).  For every $\theta\in\Theta$, the pricing PDE
\eqref{eq:rp-pricing-pde} has a unique continuous viscosity solution
in the price class $0\le u^\theta(z,t)\le K_0t$.  Moreover,
$A_0\le D^\theta\le M_D$ and, for every $(z,t)\in Q_T$,
\begin{equation}\label{eq:rp-identification}
 v^\theta(z,t)=V_0^{\theta,z,t},\qquad
 u^\theta(z,t)=\psi(V_0^{\theta,z,t})-h(y,t)-\xi(t).
\end{equation}
\end{proposition}

The well-posedness and the bounds come from \cite{Papi2002}, and the identification is the nonlinear Feynman-Kac formula of \cite{EKPQ1997} applied to the transformed equations. By \eqref{eq:rp-identification}, $\theta$ enters $u^\theta$ only through the prepayment ODE in \eqref{eq:rp-forward} and the driver $f^\theta$, and in both it appears linearly in the coefficients. Differentiating there gives $\mathcal{L}$ without differentiating $u^0$.

\subsection{Numerical illustration}\label{sec:rp-numerics}

We illustrate \eqref{eq:rp-main}, \eqref{eq:rp-ms} and \eqref{eq:rp-mean} on the
two-state example of \cite[Section~10.5]{PapiBriani2011}: $T=10$, $\mathit{MB}(0)=A_0=1$,
$\tau'=0.06$, $\tau=0.05$, $r=0.02$, $\rho=0.5$, $\sigma=0.2$, $s_0(y)=y$,
$b(x)=0.29368\tanh(0.07935-x)$, and $\eta^0(x)=0.05\,e^{-0.05/(0.04-x)}$ for $x<0.04$,
$\eta^0(x)=0$ for $x\ge0.04$. The constants are those of
\cite[Table~10.1]{PapiBriani2011}; the $\tanh$ makes $b$ bounded, and $\eta^0$ is the
smooth intensity \cite[(10.88)]{PapiBriani2011}, since the kinked one is not in
$W^{2,\infty}$. The noise is $\zeta(t)=\tfrac12\bigl(a_1+a_2\cos(2\pi t/T)\bigr)$ with
$a_1,a_2$ i.i.d., of mean $0$ and variance $1$, either symmetric ($a_i=\pm1$) or
skewed ($a_i=\sqrt2$ or $-1/\sqrt2$ with probabilities $\tfrac13$, $\tfrac23$). The two
laws give $\zeta$ the same covariance, hence the same leading term in
\eqref{eq:rp-ms}, but different third moments. Here $\varepsilon_*=1/2$ and
$1/(2\sqrt2)$, all expectations are sums over four atoms, and in
Corollary~\ref{theorem:rp-mean-variance} we take $t=T$ and $\nu=\delta_{(0.05,0)}$,
the price at issuance.

We discretize \eqref{eq:rp-pricing-pde} on $[-1,1]\times[0,2]$ by IMEX Euler, implicit
in $x$ and explicit with upwinding in $y$ ($y=0$ is an outflow boundary, $u=0$ at
$y=2$), with $321\times161$ points and $1000$ steps. $Y$ is the $\varepsilon$-derivative
of the scheme, so $u^\varepsilon-u^0-\varepsilon Y$ is the remainder of the scheme
itself. On $161\times81$ points with $500$ steps, $\mathbb E\|Y(T)\|_H^2$ changes by
$1.6\%$ and the slopes below by less than $10^{-4}$; on $[-1.5,1.5]\times[0,3]$,
$P^0(T)$ changes by $3.5\times10^{-9}$.

\begin{figure}[ht]
\centering
\includegraphics[width=\textwidth]{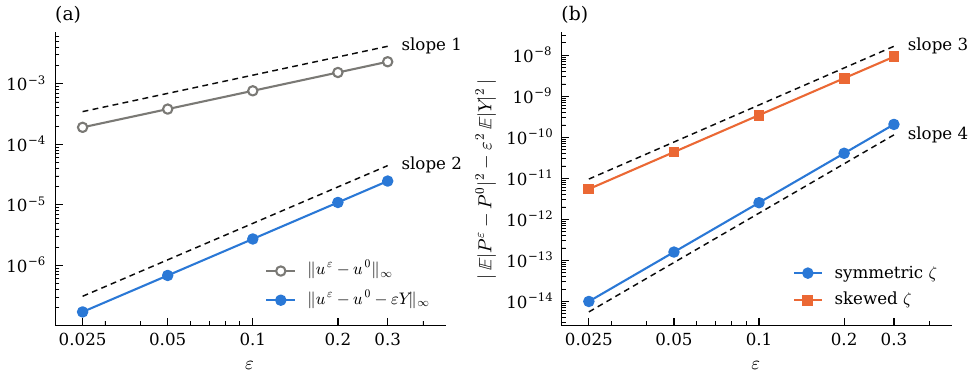}
\caption{(a) Sup norms over $[-0.2,0.3]\times[0,0.25]$ and all time steps, over the
four paths of the symmetric law. (b) Left-hand side of \eqref{eq:rp-ms} at $t=T$,
$\nu=\delta_{(0.05,0)}$. Dashed lines are reference slopes.}
\label{fig:rp}
\end{figure}

Over $0.025\le\varepsilon\le0.2$ the slopes in Figure~\ref{fig:rp}(a) are $1.00$ and
$2.00$, and at $\varepsilon=0.2$ the remainder is $0.71\%$ of $u^\varepsilon-u^0$. In
Figure~\ref{fig:rp}(b) the skewed law gives slope $2.99$, so the order
$|\varepsilon|^3$ in \eqref{eq:rp-ms} cannot be improved in general, while the
symmetric law gives $4.00$. Indeed, if $\theta\mapsto u^\theta$ is $C^3$ and
$R^\varepsilon=\varepsilon^2R_2+O(\varepsilon^3)$ with $R_2$ quadratic in $\zeta$, the
$\varepsilon^3$ term in \eqref{eq:rp-ms} is
$2\varepsilon^3\,\mathbb E\langle Y(T),R_2(T)\rangle_H$, a third moment of $\zeta$,
which is zero for the symmetric law. The variance \eqref{eq:rp-var} behaves the same
way. The mean bias is
$\mathbb EP^\varepsilon(T)-P^0(T)\approx9.4\times10^{-5}\varepsilon^2$ for both laws,
against $P^0(T)=1.1475$: centered noise in the prepayment intensity raises the mean
price at second order.

\FloatBarrier
\section{Proofs of the results}\label{sec:proofs}
\subsection{Proof of the Proposition~\ref{prop:rp-bsde-identification}}
First, we use Theorem~4.2 and Corollary~4.4 in \cite{Papi2002} to show that 
\eqref{eq:rp-pricing-pde} has a unique continuous viscosity solution. Here we verify the
conditions (P1)--(P3), (4.35) and (4.36) of \cite{Papi2002}. By (A1) and
$\|\theta\|_\infty\le\frac12$, it is easy to see that
\[
 \|\eta^\theta\|_\infty\le\tfrac32\|\eta^0\|_\infty,\qquad
 \|\nabla\eta^\theta(\cdot,t)\|_\infty\le\tfrac32\sup_{s\in[0,T]}\|\nabla\eta^0(\cdot,s)\|_\infty .
\]
Since $\mu^\theta=(b,\eta^\theta)$ and $b$ is bounded and Lipschitz in $x$, (P1) holds.
Recall \eqref{eq:rp-balance-and-price} and compute the derivatives: 
\[
h_y=-s_0'(y)h,\qquad h_t=-\mathit{MB}'(T-t)e^{-s_0(y)}, \qquad 
h_x=h_{xx}=0
\]
Hence $\operatorname{tr}(\Sigma\Sigma^\top D_z^2h)=\sigma^2(t)h_{xx}=0$, where $\Sigma=(\sigma,0)^\top$. By (A1), $h\in C^{2,1}$ is nonnegative, and $h$, $\nabla h$, $h_t$ are 
bounded and spatially Lipschitz. Since also $\xi\in C^1$ and $\xi\ge A_0>0$, (P2) holds. (P3) is trivial since $U_0\equiv0$. It is easy to conclude 
\[
|h(y,t)-h(y,s)|\le\|\mathit{MB'}\|_\infty e^{-\inf s_0}|t-s|,
\]
which is (4.35). Since $\tau\ge r$, $h\ge0$ and $h(y,t)\to0$ as $y\to\infty$, we have 
$\inf_z(\tau-r)h(z,t)=0$ for each $t$, so the barriers in \cite[(4.34)]{Papi2002} are
$\underline k=0$ and $\overline k(t)=K_0t$. Thus $\xi+h+\underline k\ge A_0>0$, and (4.36)
holds.

Applying Theorem~4.2 and Corollary~4.4 gives a unique continuous viscosity solution
$u^\theta$ on $\mathbb R^2\times[0,T)$ between the barriers, that is,
$0\le u^\theta(z,t)\le K_0t$. Since $D^\theta=u^\theta+h+\xi$,
$\xi(t)=A_0e^{r(T-t)}$ and $M_D=K_0T+\|h\|_\infty+\|\xi\|_\infty$, we conclude
$A_0\le D^\theta\le M_D$.

Next we use the power transform to cancel the nonlinear term in \eqref{eq:rp-pricing-pde}
By definition $D^\theta=u^\theta +h+\xi$, so
\[
 D_t^\theta=u_t^\theta+h_t+\xi',\quad D_x^\theta=u_x^\theta,\quad D_{xx}^\theta=u_{xx}^\theta,\quad
 u_y^\theta=D_y^\theta-h_y,\quad u^\theta=D^\theta-h-\xi .
\]
Plugging into \eqref{eq:rp-pricing-pde} and using $\xi'=-r\xi$ gives
\begin{equation}\label{eq:rp-D-pde}
 D_t^\theta=aD_{xx}^\theta+bD_x^\theta+\eta^\theta D_y^\theta-k\frac{(D_x^\theta)^2}{D^\theta}
 +\tau h+h_t-\eta^\theta h_y-rD^\theta .
\end{equation}
Recall from \eqref{eq:rp-transform} that $\Phi'(d)=d^{-2\rho}$,
$\Phi''(d)=-2\rho d^{-2\rho-1}$ and $v^\theta=\Phi(D^\theta)$. Computing the derivatives gives
\[
 v_t^\theta=\Phi'(D^\theta)D_t^\theta,\quad v_x^\theta=\Phi'(D^\theta)D_x^\theta,\quad
 v_y^\theta=\Phi'(D^\theta)D_y^\theta,\quad
 v_{xx}^\theta=\Phi''(D^\theta)(D_x^\theta)^2+\Phi'(D^\theta)D_{xx}^\theta .
\]
Plugging into \eqref{eq:rp-D-pde} gives
\[
 v_t^\theta-av_{xx}^\theta-bv_x^\theta-\eta^\theta v_y^\theta
 =\Bigl(-\frac{k}{D^\theta}\Phi'(D^\theta)-a\Phi''(D^\theta)\Bigr)(D_x^\theta)^2
 +\Phi'(D^\theta)\bigl(\tau h+h_t-\eta^\theta h_y-rD^\theta\bigr).
\]
Note that
\[
 -\frac{k}{D^\theta}\Phi'(D^\theta)-a\Phi''(D^\theta)
 =-\frac{k}{D^\theta}(D^\theta)^{-2\rho}+2\rho a(D^\theta)^{-2\rho-1}
 =(2\rho a-k)(D^\theta)^{-2\rho-1}=0
\]
since $k=2\rho a$. Rearranging,
\begin{equation}\label{eq:rp-v-pde-raw}
 v_t^\theta=av_{xx}^\theta+bv_x^\theta+\eta^\theta v_y^\theta
 +(D^\theta)^{-2\rho}\bigl(\tau h+h_t-\eta^\theta h_y-rD^\theta\bigr).
\end{equation}
Since $A_0\le D^\theta\le M_D$, $v^\theta$ takes values in $I$, on which $\psi=\Phi^{-1}$, so
$D^\theta=\psi(v^\theta)$. By \eqref{eq:rp-transformed-driver}, rearranging
\eqref{eq:rp-v-pde-raw} we conclude
\begin{equation}\label{eq:rp-v-pde}
\begin{cases}
 v_t^\theta=av_{xx}^\theta+bv_x^\theta+\eta^\theta v_y^\theta+f^\theta(z,t,v^\theta),\\
 v^\theta(z,0)=G(y).
\end{cases}
\end{equation}
The computation above is for smooth functions. For $u^\theta$ it holds in the viscoity sense, if a smooth $\varphi$ touches $v^\theta$ from above(below) at a ppoint, then 
$\Phi^{-1}\circ\varphi-h-\xi$ touches $u^\theta$ from above (below) there, and the two viscosity inequalities differ by the factor $\Phi'(D^\theta)>0$. So $v^\theta $ is a bounded continuous viscosity solution of \eqref{eq:rp-v-pde} on $\mathbb R^2\times[0,T)$.
Recall \eqref{eq:rp-forward}. By (A1), $b$ and $\eta^\theta$ are bounded and globally
Lipschitz, and $\sigma$ is continuous and bounded. By the existence and uniqueness theorem
for SDEs, \eqref{eq:rp-forward} has a unique solution.

Finally, we apply Theorem~2.1 and Theorem~4.2 in \cite{EKPQ1997} to show that the value
function of the BSDE corresponds to the viscosity solution of \eqref{eq:rp-v-pde}. Since
$G(y)=\Phi(h(y,0)+\xi(0))$,
\[
 G_y(y)=\bigl(h(y,0)+\xi(0)\bigr)^{-2\rho}h_y(y,0).
\]
Since $h+\xi\ge A_0>0$, we have $|G_y(y)|\le A_0^{-2\rho}\|h_y\|_\infty<\infty$, so $G$ is
Lipschitz. $G$ is also bounded, since $h(\cdot,0)+\xi(0)$ takes values in
$[A_0,\|h\|_\infty+\xi(0)]$. Since $\psi$ is bounded and bounded away from zero, $\psi'$ is
bounded, and $h$, $h_t$, $h_y$, $\eta^\theta$ are bounded and spatially Lipschitz, $f^\theta$
is bounded and
\[
 |f^\theta(z,t,w)-f^\theta(z',t,w')|\le L\bigl(|z-z'|+|w-w'|\bigr).
\]
By \cite[Theorem~2.1]{EKPQ1997}, \eqref{eq:rp-bsde} admits a unique square-integrable
solution $(V^\theta,Q^\theta)$.

Note that the PDE in \cite[Theorem~4.2]{EKPQ1997} is a terminal-value PDE, so we need to
reverse time. Let $\widehat v(s,z)=v^\theta(z,T-s)$, and similarly
\[
 \widehat\mu(s,z)=\begin{pmatrix}b(x,T-s)\\ \eta^\theta(z,T-s)\end{pmatrix},\qquad
 \widehat\Sigma(s)=\Sigma(T-s),\qquad
 \widehat f(s,z,w,q)=f^\theta(z,T-s,w),
\]
with terminal value $g(z)=G(y)$. The forward--backward system in \cite{EKPQ1997} with these
data, started at time $T-t$ from $z$, is \eqref{eq:rp-forward}--\eqref{eq:rp-bsde} after
shifting the time origin by $T-t$, and $Y^{s,z}_s$ is deterministic
\cite[Proposition~4.2]{EKPQ1997}. So the function $u(s,z)=Y^{s,z}_s$ in
\cite[Theorem~4.2]{EKPQ1997} is $V_0^{\theta,z,T-s}$. By the same estimate,
\[
 |\widehat f(s,z,w,q)-\widehat f(s,z',w,q)|\le L|z-z'| .
\]
The BSDE is scalar and $g$, $\widehat f$ are Lipschitz in $z$. Choose $m_R(\delta)=L\delta$
in \cite[(4.12)]{EKPQ1997}; then Theorem~4.2 gives that $(z,t)\mapsto V_0^{\theta,z,t}$ is
continuous on $Q_T$ and is the unique viscosity solution of \eqref{eq:rp-v-pde}.

By the power transform, $v^\theta$ is a bounded continuous viscosity solution of
\eqref{eq:rp-v-pde} with the same initial value. Since $u^\theta$ is so far only defined for
$t<T$, we apply uniqueness on $\mathbb R^2\times[0,T_0]$ for each $T_0<T$, and get
$v^\theta(z,t)=V_0^{\theta,z,t}$ for $t<T$. The value function is continuous on $Q_T$ and
takes values in $I$, so $\psi(V_0^{\theta,z,t})-h-\xi$ extends $u^\theta$ continuously to
$t=T$; by continuity the bounds and the uniqueness in the price class also hold on $Q_T$.
Since $v^\theta=\Phi(D^\theta)$, rearranging gives
\[
 u^\theta(z,t)=\psi(V_0^{\theta,z,t})-h(y,t)-\xi(t).
\]
\subsection{Proof of the Theorem~\ref{theorem:rp-first-order}}
By Proposition~\ref{prop:rp-bsde-identification}, $u^\theta=\psi(V^\theta_0)-h-\xi$, so the expansion of $u^\theta$ follows from the expansions of the forward variation $\mathcal A^\theta-\mathcal A^0$ and of the backward variation $V^\theta-V^0$, transferred to $u^\theta$ through $\psi$.
Fix $\theta\in\Theta$ and $(z,t)\in Q_T$, and write $\delta=\|\theta\|_\infty$.
Throughout, $C$ denotes a constant that depends only on (A1) and may change from line to line; it does not depend on $s$, $\omega$, $(z,t)$ or $\theta$.
We write $O(\delta^2)$ for a term bounded in absolute value by $C\delta^2$.

We first show that the forward state admits a first-order expansion in $\theta$, uniformly in time, in the initial state and in $\omega$.
The proof only uses the ODE for $\mathcal A^\theta$, a Taylor expansion of $\eta^0$ and Gronwall's inequality.
Combining \eqref{eq:re-scenario-set} and \eqref{eq:rp-forward} gives
\[
\frac{d\mathcal A^\theta_s}{ds}=\eta^0(X_s,\mathcal A^\theta_s,t-s)\bigl[1+\theta(t-s)\bigr].
\]
Define $\Delta\mathcal A_s=\mathcal A^\theta_s-\mathcal A^0_s$.
Since $X$ does not depend on $\theta$ and both processes start at $y$, we have $\Delta\mathcal A_0=0$ and
\begin{align*}
\frac{d}{ds}\Delta\mathcal A_s
&=\eta^0(X_s,\mathcal A^\theta_s,t-s)\bigl[1+\theta(t-s)\bigr]-\eta^0(X_s,\mathcal A^0_s,t-s)\\
&=\eta^0(X_s,\mathcal A^\theta_s,t-s)-\eta^0(X_s,\mathcal A^0_s,t-s)+\eta^0(X_s,\mathcal A^\theta_s,t-s)\,\theta(t-s).
\end{align*}
By (A1), $\eta^0$ is bounded by a constant $M_\eta$ and Lipschitz in $y$ with a constant $L_\eta$, uniformly in $(x,t)$.
Taking absolute values gives
\begin{align*}
\Bigl|\frac{d}{ds}\Delta\mathcal A_s\Bigr|
&\le\bigl|\eta^0(X_s,\mathcal A^\theta_s,t-s)-\eta^0(X_s,\mathcal A^0_s,t-s)\bigr|+\bigl|\eta^0(X_s,\mathcal A^\theta_s,t-s)\,\theta(t-s)\bigr|\\
&\le L_\eta|\Delta\mathcal A_s|+M_\eta|\theta(t-s)|
\le L_\eta|\Delta\mathcal A_s|+M_\eta\delta .
\end{align*}
Note that $\Delta\mathcal A_s=\int_0^s\frac{d}{dr}\Delta\mathcal A_r\,dr$, so
\[
|\Delta\mathcal A_s|\le\int_0^s\Bigl|\frac{d}{dr}\Delta\mathcal A_r\Bigr|\,dr
\le\int_0^s\bigl(L_\eta|\Delta\mathcal A_r|+M_\eta\delta\bigr)\,dr
=L_\eta\int_0^s|\Delta\mathcal A_r|\,dr+M_\eta\delta s .
\]
Gronwall's inequality gives $|\Delta\mathcal A_s|\le M_\eta\delta s\,e^{L_\eta s}\le M_\eta Te^{L_\eta T}\delta$.
Therefore
\begin{equation}\label{eq:rp-F1}
\sup_{0\le s\le t}|\Delta\mathcal A_s|\le C\delta .
\end{equation}

Next we introduce the first-order variation process $J_s[\theta]$, defined by
\[
\frac{d}{ds}J_s[\theta]=\eta^0_y(X_s,\mathcal A^0_s,t-s)\,J_s[\theta]+\eta^0(X_s,\mathcal A^0_s,t-s)\,\theta(t-s),\qquad J_0[\theta]=0 .
\]
Since $|\eta^0_y|\le L_\eta$ and $|\eta^0|\le M_\eta$, the same argument as for $\Delta\mathcal A$ gives $\sup_{0\le s\le t}|J_s[\theta]|\le C_J\delta$.

Next we define $E^A_s:=\Delta\mathcal A_s-J_s[\theta]$. Then
\begin{equation}\label{eq:rp-EA}
\begin{aligned}
\frac{d}{ds}E^A_s
&=\frac{d}{ds}\Delta\mathcal A_s-\frac{d}{ds}J_s[\theta]\\
&=\eta^0(X_s,\mathcal A^\theta_s,t-s)-\eta^0(X_s,\mathcal A^0_s,t-s)-\eta^0_y(X_s,\mathcal A^0_s,t-s)\,J_s[\theta]\\
&\quad+\bigl[\eta^0(X_s,\mathcal A^\theta_s,t-s)-\eta^0(X_s,\mathcal A^0_s,t-s)\bigr]\theta(t-s).
\end{aligned}
\end{equation}
By Taylor's theorem,
\[
\eta^0(X_s,\mathcal A^\theta_s,t-s)-\eta^0(X_s,\mathcal A^0_s,t-s)=\eta^0_y(X_s,\mathcal A^0_s,t-s)\,\Delta\mathcal A_s+\mathcal T^A_s ,
\]
and by (A1) the remainder satisfies
\[
|\mathcal T^A_s|\le\tfrac12\|\eta^0_{yy}\|_\infty|\Delta\mathcal A_s|^2\le\tfrac12\|\eta^0_{yy}\|_\infty C^2\delta^2 .
\]
Hence $\mathcal T^A_s=O(\delta^2)$.
Similarly, by \eqref{eq:rp-F1}, the last term in \eqref{eq:rp-EA} satisfies
\[
\bigl|\bigl[\eta^0(X_s,\mathcal A^\theta_s,t-s)-\eta^0(X_s,\mathcal A^0_s,t-s)\bigr]\theta(t-s)\bigr|\le L_\eta|\Delta\mathcal A_s|\,\delta\le C\delta^2 .
\]
Therefore
\[
\frac{d}{ds}E^A_s=\eta^0_y(X_s,\mathcal A^0_s,t-s)\,E^A_s+O(\delta^2).
\]
Since $E^A_0=0$ and $|\eta^0_y|\le L_\eta$, Gronwall's inequality gives $\sup_{0\le s\le t}|E^A_s|\le C\delta^2$.
Combining the results above gives
\begin{equation}\label{eq:rp-F2}
\sup_{0\le s\le t}\bigl|\mathcal A^\theta_s-\mathcal A^0_s-J_s[\theta]\bigr|\le C\delta^2 .
\end{equation}

Next we show the expansion of the backward variation.
Define $\Delta V_s=V^\theta_s-V^0_s$ and $\Delta Q_s=Q^\theta_s-Q^0_s$.
Following \eqref{eq:rp-bsde},
\begin{equation}\label{eq:rp-B1}
\begin{aligned}
\Delta V_s={}&G(\mathcal A^\theta_t)-G(\mathcal A^0_t)
+\int_s^t\bigl[f^\theta(S^\theta_r,t-r,V^\theta_r)-f^0(S^0_r,t-r,V^0_r)\bigr]\,dr\\
&-\int_s^t\Delta Q_r\,dB_r ,
\end{aligned}
\end{equation}
where $S^\theta_r=(X_r,\mathcal A^\theta_r)$ and $S^0_r=(X_r,\mathcal A^0_r)$ are as in \eqref{eq:rp-forward}.
Recall from the proof of Proposition~\ref{prop:rp-bsde-identification} that $G$ is globally Lipschitz and that $f^\theta$ is Lipschitz in $(z,w)$ with a constant $L$ that does not depend on $\theta\in\Theta$.
Write $\Delta G=G(\mathcal A^\theta_t)-G(\mathcal A^0_t)$.
Then $|\Delta G|\le L_G|\mathcal A^\theta_t-\mathcal A^0_t|$, and hence, by \eqref{eq:rp-F1}, $|\Delta G|\le C_G\delta$.

We split the difference of the drivers as
\begin{align*}
f^\theta(S^\theta_r,t-r,V^\theta_r)-f^0(S^0_r,t-r,V^0_r)
&=\underbrace{f^\theta(S^\theta_r,t-r,V^\theta_r)-f^\theta(S^\theta_r,t-r,V^0_r)}_{=:\beta_r}\\
&\quad+\underbrace{f^\theta(S^\theta_r,t-r,V^0_r)-f^0(S^0_r,t-r,V^0_r)}_{=:\lambda_r}.
\end{align*}
Expand $\lambda_r$ to
\[
\lambda_r=\underbrace{f^\theta(S^\theta_r,t-r,V^0_r)-f^\theta(S^0_r,t-r,V^0_r)}_{=:\mathrm{I}_r}
+\underbrace{f^\theta(S^0_r,t-r,V^0_r)-f^0(S^0_r,t-r,V^0_r)}_{=:\mathrm{II}_r},
\]
so that $|\lambda_r|\le|\mathrm{I}_r|+|\mathrm{II}_r|\le C\delta+C\delta\le C_\lambda\delta$.
Indeed, since $f^\theta$ is spatially Lipschitz, $|\mathrm{I}_r|\le L|\Delta\mathcal A_r|\le C\delta$ by \eqref{eq:rp-F1}.
Plugging \eqref{eq:re-scenario-set} into \eqref{eq:rp-transformed-driver}, it is easy to see that
\[
\mathrm{II}_r=-\psi(V^0_r)^{-2\rho}\,\eta^0(S^0_r,t-r)\,h_y(\mathcal A^0_r,t-r)\,\theta(t-r),
\]
so $|\mathrm{II}_r|\le C\delta$, since $\psi$ is bounded away from zero and $\eta^0$, $h_y$ are bounded.

Now we deal with $\beta_r$. We define
\[
a_r=
\begin{cases}
\dfrac{\beta_r}{\Delta V_r}, & \Delta V_r\neq0,\\[1.5ex]
0, & \Delta V_r=0 .
\end{cases}
\]
Since $|\beta_r|\le L|\Delta V_r|$, we have $|a_r|\le L$, and $a$ is progressively measurable.
So $\beta_r=a_r\Delta V_r$.
Combining the results above with \eqref{eq:rp-B1}, we get
\begin{equation}\label{eq:rp-B2}
\Delta V_s=\Delta G+\int_s^t\bigl(a_r\Delta V_r+\lambda_r\bigr)\,dr-\int_s^t\Delta Q_r\,dB_r .
\end{equation}
Applying an integrating factor and conditional expectation gives
\[
\Delta V_s=\widehat{\mathbb E}\Bigl[e^{\int_s^ta_\ell\,d\ell}\,\Delta G+\int_s^te^{\int_s^ra_\ell\,d\ell}\,\lambda_r\,dr\ \Big|\ \mathcal F_s\Bigr],
\]
where $(\mathcal F_s)$ is the filtration of $B$.
Indeed, by It\^o's formula and \eqref{eq:rp-B2},
\[
d\bigl(e^{\int_s^ra_\ell\,d\ell}\Delta V_r\bigr)=e^{\int_s^ra_\ell\,d\ell}\bigl(-\lambda_r\,dr+\Delta Q_r\,dB_r\bigr),
\]
and the stochastic integral has zero conditional expectation since $a$ is bounded and $\Delta Q$ is square integrable.
Since $|a_\ell|\le L$, we have $0<e^{\int_s^ra_\ell\,d\ell}\le e^{L(r-s)}\le e^{LT}$, so
\[
|\Delta V_s|\le\widehat{\mathbb E}\Bigl[e^{LT}C_G\delta+\int_s^te^{LT}C_\lambda\delta\,dr\ \Big|\ \mathcal F_s\Bigr]
=e^{LT}\bigl[C_G+C_\lambda(t-s)\bigr]\delta\le e^{LT}\bigl(C_G+TC_\lambda\bigr)\delta .
\]
Define $C'=e^{LT}(C_G+TC_\lambda)$. This gives
\begin{equation}\label{eq:rp-dV}
\sup_{0\le s\le t}|\Delta V_s|\le C'\delta\qquad\text{a.s.}
\end{equation}

Next, define
\[
c_r=\partial_wf^0(S^0_r,t-r,V^0_r),\qquad
d_r=\partial_yf^0(S^0_r,t-r,V^0_r),\qquad
e_r=-\psi(V^0_r)^{-2\rho}\,\eta^0(S^0_r,t-r)\,h_y(\mathcal A^0_r,t-r),
\]
where $\partial_y$ is taken with $w$ fixed.
Note that $|c_r|+|d_r|+|e_r|\le C$; in particular $|c_r|\le C_c$ for a constant $C_c$.
Now define $R_s[\theta]$ by
\begin{equation}\label{eq:rp-B3}
R_s[\theta]=G_y(\mathcal A^0_t)J_t[\theta]
+\int_s^t\bigl(c_rR_r[\theta]+d_rJ_r[\theta]+e_r\theta(t-r)\bigr)\,dr
-\int_s^t\widetilde Q_r[\theta]\,dB_r .
\end{equation}
This is a linear BSDE with bounded coefficients, so it has a unique solution $(R[\theta],\widetilde Q[\theta])$ by \cite[Theorem~2.1]{EKPQ1997}.
Note that $|G_y(\mathcal A^0_t)J_t[\theta]|\le\|G_y\|_\infty|J_t[\theta]|\le\|G_y\|_\infty C_J\delta$ and
\[
|d_rJ_r[\theta]+e_r\theta(t-r)|\le|d_r||J_r[\theta]|+|e_r||\theta(t-r)|\le\bigl(\|d\|_\infty C_J+\|e\|_\infty\bigr)\delta .
\]
By the same integrating factor argument as for $\Delta V$, now with $|c_r|\le C_c$,
\[
\sup_{0\le s\le t}|R_s[\theta]|\le e^{C_cT}\bigl(\|G_y\|_\infty C_J+T\|d\|_\infty C_J+T\|e\|_\infty\bigr)\delta .
\]
Define $C_R=e^{C_cT}(\|G_y\|_\infty C_J+T\|d\|_\infty C_J+T\|e\|_\infty)$.
We conclude $\sup_{0\le s\le t}|R_s[\theta]|\le C_R\delta$ a.s.

Now let
\begin{align*}
\mathcal T^V_r:={}&f^\theta(S^\theta_r,t-r,V^\theta_r)-f^0(S^0_r,t-r,V^0_r)\\
&-c_r\Delta V_r-d_r\Delta\mathcal A_r-e_r\theta(t-r).
\end{align*}
By \eqref{eq:re-scenario-set} and \eqref{eq:rp-transformed-driver},
\[
f^\theta(z,t,w)=f^0(z,t,w)+\theta(t)\,g(z,t,w),\qquad g(z,t,w)=-\psi(w)^{-2\rho}\eta^0(z,t)h_y(y,t),
\]
and $e_r=g(S^0_r,t-r,V^0_r)$.
Since $S^\theta_r$ and $S^0_r$ have the same first component $X_r$, we apply Taylor's theorem to $f^0$ in $(y,w)$.
By (A1) and the bounds on $\psi,\psi',\psi''$, the first derivatives of $f^0$ in $(y,w)$ and the function $g$ are Lipschitz in $(y,w)$, uniformly in $(x,t)$.
Hence, by \eqref{eq:rp-F1} and \eqref{eq:rp-dV},
\[
|\mathcal T^V_r|\le C\bigl[|\Delta\mathcal A_r|^2+|\Delta V_r|^2+|\theta(t-r)|\bigl(|\Delta\mathcal A_r|+|\Delta V_r|\bigr)\bigr]\le C\delta^2 .
\]
So we obtain
\[
f^\theta(S^\theta_r,t-r,V^\theta_r)-f^0(S^0_r,t-r,V^0_r)=c_r\Delta V_r+d_r\Delta\mathcal A_r+e_r\theta(t-r)+\mathcal T^V_r .
\]
Observe that, since $G_y$ is Lipschitz by (A1),
\begin{align*}
H_t:={}&G(\mathcal A^\theta_t)-G(\mathcal A^0_t)-G_y(\mathcal A^0_t)J_t[\theta]\\
={}&G_y(\mathcal A^0_t)\Delta\mathcal A_t+O(\delta^2)-G_y(\mathcal A^0_t)J_t[\theta]\\
={}&G_y(\mathcal A^0_t)\bigl(\Delta\mathcal A_t-J_t[\theta]\bigr)+O(\delta^2)\\
={}&O(\delta^2)+O(\delta^2),
\end{align*}
By Taylor's theorem, \eqref{eq:rp-F1} and \eqref{eq:rp-F2}, $|H_t|\le C\delta^2$.

Now subtract the two BSDEs. Let
\begin{equation}\label{eq:rp-star}
E^V_s:=\Delta V_s-R_s[\theta].
\end{equation}
Subtracting \eqref{eq:rp-B3} from \eqref{eq:rp-B1} and using the expansion of the driver above and $\Delta\mathcal A_r-J_r[\theta]=E^A_r$, we obtain
\[
E^V_s=H_t+\int_s^t\bigl(c_rE^V_r+d_rE^A_r+\mathcal T^V_r\bigr)\,dr-\int_s^t\bigl(\Delta Q_r-\widetilde Q_r[\theta]\bigr)\,dB_r .
\]
By the results above, $|H_t|\le C\delta^2$ and, by \eqref{eq:rp-F2}, $|d_rE^A_r+\mathcal T^V_r|\le C\delta^2$.
By the same integrating factor argument as for $\Delta V$, it is easy to conclude
\[
\sup_{0\le s\le t}|E^V_s|\le C\delta^2\qquad\text{a.s.}
\]
Plugging back into \eqref{eq:rp-star} gives
\begin{equation}\label{eq:rp-B4}
\sup_{0\le s\le t}\bigl|V^\theta_s-V^0_s-R_s[\theta]\bigr|\le C\delta^2\qquad\text{a.s.}
\end{equation}

Finally, we transfer the expansion to the price.
Applying Proposition~\ref{prop:rp-bsde-identification}, by \eqref{eq:rp-identification},
\[
u^\theta(z,t)=\psi\bigl(V^{\theta,z,t}_0\bigr)-h(y,t)-\xi(t),
\qquad\text{so}\qquad
u^\theta(z,t)-u^0(z,t)=\psi(V^\theta_0)-\psi(V^0_0).
\]
Applying a Taylor expansion to $\psi$, whose second derivative is bounded, and using \eqref{eq:rp-dV} and \eqref{eq:rp-B4} gives
\[
u^\theta(z,t)-u^0(z,t)=\psi'(V^0_0)\,\Delta V_0+O(\delta^2)=\psi'(V^0_0)\,R_0[\theta]+O(\delta^2).
\]
The same variational equations define $J[\theta]$ and $R[\theta]$
for every $\theta\in C([0,T])$.
Now define the price sensitivity operator
\[
(\mathcal L\theta)(z,t):=\psi'\bigl(V^{0,z,t}_0\bigr)\,R^{z,t}_0[\theta].
\]
Since the constants above do not depend on $(z,t)$, we get $\|u^\theta-u^0-\mathcal L\theta\|_\infty\le C_2\|\theta\|_\infty^2$, which is \eqref{eq:rp-Frechet}.
Since $J[\theta]$ and $R[\theta]$ are linear in $\theta$ and $\sup_s|R_s[\theta]|\le C_R\|\theta\|_\infty$, the operator $\mathcal L$ is linear and $\|\mathcal L\theta\|_\infty\le\|\psi'\|_\infty C_R\|\theta\|_\infty$.
Moreover, \eqref{eq:rp-Frechet} applied to $\theta/n$ gives $\|n(u^{\theta/n}-u^0)-\mathcal L\theta\|_\infty\le C_2\|\theta\|_\infty^2/n$, so $\mathcal L\theta$ is a uniform limit of bounded continuous functions and $\mathcal L\theta\in C_b(Q_T)$.

Take $\theta=\varepsilon\zeta(\omega)$ and $Y=\mathcal L\zeta$.
For $|\varepsilon|\le\varepsilon_*$, $\varepsilon\zeta(\omega)\in\Theta$ for a.e.\ $\omega$ and $\mathcal L(\varepsilon\zeta)=\varepsilon Y$, so
\[
\|u^\varepsilon-u^0-\varepsilon Y\|_\infty\le C_2\varepsilon^2N(\omega)^2 ,
\]

which is \eqref{eq:rp-main}.
Since $\mathcal L$ is bounded and linear and $C([0,T])$ is separable, $Y$ is strongly measurable as an $\mathcal X$-valued random variable.
The same holds for $u^\varepsilon$, repeating the proofs of \eqref{eq:rp-F1} and \eqref{eq:rp-dV} with any $\theta_2\in\Theta$ in place of $0$ gives $\|u^{\theta_1}-u^{\theta_2}\|_\infty\le C\|\theta_1-\theta_2\|_\infty$ for $\theta_1,\theta_2\in\Theta$, since $\psi$ is Lipschitz.
Using $\|v\|_{\mathcal X}\le\|v\|_\infty$ from \eqref{eq:rp-spaces} and taking the expectation of the square,
\[
\|u^\varepsilon-u^0-\varepsilon Y\|_{L^2(\Omega;\mathcal X)}\le C_2\varepsilon^2\bigl(\mathbb EN^4\bigr)^{1/2},
\]
which is \eqref{eq:rp-main-L2}. This completes the proof.
\subsection{Proof of the Corollary~\ref{theorem:rp-mean-variance}}
Fix $t\in[0,T]$ and $|\varepsilon|\le\varepsilon_*$.
Since $P^\varepsilon-P^0=u^\varepsilon-u^0$, Theorem~\ref{theorem:rp-first-order} gives
\begin{equation}\label{eq:rp-cor-split}
P^\varepsilon(t)-P^0(t)=\varepsilon Y(t)+R^\varepsilon(t),\qquad Y=\mathcal L\zeta .
\end{equation}
Let $M=\|\mathcal L\|_{C([0,T])\to C_b(Q_T)}$.
Since $\|\zeta\|_\infty\le K$, \eqref{eq:rp-main} gives, almost surely,
\begin{equation}\label{eq:rp-cor-bounds}
\|Y(t)\|_H\le\|Y\|_\infty\le MK,\qquad
\|R^\varepsilon(t)\|_H\le\|R^\varepsilon\|_\infty\le C_2K^2\varepsilon^2 .
\end{equation}
By Theorem~\ref{theorem:rp-first-order}, $Y(t)$ and $R^\varepsilon(t)$ are strongly measurable, and by \eqref{eq:rp-cor-bounds} they are bounded, so all expectations below are Bochner integrals.

Expanding the square in \eqref{eq:rp-cor-split} and taking expectations,
\[
\mathbb E\|P^\varepsilon(t)-P^0(t)\|_H^2-\varepsilon^2\mathbb E\|Y(t)\|_H^2
=2\varepsilon\,\mathbb E\langle Y(t),R^\varepsilon(t)\rangle_H+\mathbb E\|R^\varepsilon(t)\|_H^2 .
\]
By \eqref{eq:rp-cor-bounds}, the right-hand side is bounded in absolute value by $2MC_2K^3|\varepsilon|^3+C_2^2K^4\varepsilon^4$.
Since $|\varepsilon|\le\varepsilon_*=(2K)^{-1}$, we have $K^4\varepsilon^4\le\frac12K^3|\varepsilon|^3$, so \eqref{eq:rp-ms} holds with
\[
C_3=K^3\bigl(2MC_2+\tfrac12C_2^2\bigr).
\]

Since $P^0$ is deterministic, taking expectations in \eqref{eq:rp-cor-split} gives $\mathbb EP^\varepsilon(t)-P^0(t)-\varepsilon\mathbb EY(t)=\mathbb ER^\varepsilon(t)$, so
\[
\|\mathbb EP^\varepsilon(t)-P^0(t)-\varepsilon\mathbb EY(t)\|_H\le\mathbb E\|R^\varepsilon(t)\|_H\le C_2K^2\varepsilon^2 ,
\]
which is \eqref{eq:rp-mean}.

Finally, subtracting this identity from \eqref{eq:rp-cor-split} gives $P^\varepsilon(t)-\mathbb EP^\varepsilon(t)=\varepsilon\widetilde Y+\widetilde R^\varepsilon$, where $\widetilde Y=Y(t)-\mathbb EY(t)$ and $\widetilde R^\varepsilon=R^\varepsilon(t)-\mathbb ER^\varepsilon(t)$.
Expanding the square as before and using the definition of $\operatorname{Var}_H$,
\[
\operatorname{Var}_H(P^\varepsilon(t))-\varepsilon^2\operatorname{Var}_H(Y(t))
=2\varepsilon\,\mathbb E\langle\widetilde Y,\widetilde R^\varepsilon\rangle_H+\mathbb E\|\widetilde R^\varepsilon\|_H^2 .
\]
Since $\mathbb E\|V-\mathbb EV\|_H^2\le\mathbb E\|V\|_H^2$, \eqref{eq:rp-cor-bounds} gives $\mathbb E\|\widetilde Y\|_H^2\le M^2K^2$ and $\mathbb E\|\widetilde R^\varepsilon\|_H^2\le C_2^2K^4\varepsilon^4$, and the Cauchy--Schwarz inequality in $L^2(\Omega;H)$ bounds the right-hand side by $C_3|\varepsilon|^3$, as in the proof of \eqref{eq:rp-ms}. This proves \eqref{eq:rp-var}.

\end{document}